\documentclass[12pt]{amsart}
\usepackage[utf8]{inputenc}
\usepackage[english]{babel}
\addto\captionsenglish{}
\usepackage[margin=2cm]{caption}
\usepackage{fancyhdr}
\usepackage[colorlinks=true,unicode=true]{hyperref}
\usepackage{MnSymbol}
\usepackage{wrapfig}
\usepackage{tikz}
\usepackage{pgfplots}
\pgfplotsset{compat=1.17}
\usepackage[backend=biber,
maxnames=10,
doi=false,
url=false,
style=alphabetic,
isbn=false]{biblatex}
\input{defs.tex}

\makeatletter
\let\c@equation\c@subsection

\makeatother
\newtheorem{theorem}[subsection]{Theorem}
\newtheorem{lemma}[subsection]{Lemma}
\newtheorem{proposition}[subsection]{Proposition}

\title{Reverse Bernstein Inequality on the Circle}
\author[P. Joharinad]{Parvaneh Joharinad}
\author[J. Jost]{Jürgen Jost}
\author[S. Lim]{Sunhyuk Lim}
\author[R. Matveev]{Rostislav Matveev}
\address{Max-Plank Institut für Mathematik in den
  Naturwissenschaften\\
  Inselstraße 22, 04103 Leipzig, Germany}

\begin{document}
\begin{abstract}
  The more then hundred years old Bernstein inequality states that the
  $L^{\infty}$-norm of the derivative of a trigonometric polynomial of
  fixed degree can be bounded from above by $L^{\infty}$-norm of the
  polynomial itself. The {\em reversed Bernstein inequality}, that we
  prove in this note, says that the reverse inequality holds for
  functions in the orthogonal complement of the space of polynomials
  of fixed degree.

  In fact, we derived a more general result for the lower bounds on
  higher derivatives. These bounds are better then those obtained by
  applying bound for the first derivative successively several times.

  Rostislav Matveev is very grateful to Evgeniy Abakumov for fruitful
  discussions about this topic.
\end{abstract}

\maketitle

\section{Introduction and the statement}
We have to come to consider reverse Bernstein inequality while
thinking about homotopy types of neighborhoods of the circle in its
hyperconvex hull.  While at this project is ongoing, we decided
to write this small note about the inequality itself without
applications.
\medskip

Define $S^{1}:=\Rbb/2\pi\Zbb$ with the angular metric
\[
  d(x,y):=d_{\Rbb}(x,y+2\pi\Zbb)
\]
and the normalized measure
$\d\mu:=\frac{1}{2\pi}\d x$. We will use angular coordinates on the
circle in the range $(-\pi,\pi]$. Let $L^{2}:=L^{2}(S^{1},\mu)$ be the
space of complex-valued $L^{2}$-functions on $S^{1}$.

For a number $k\in \Nbb$ define the space of trigonometric
polynomials and its orthogonal complement by
\begin{align*}
  \Pcal_{k}
  &:=
    \left\llangle \ebf^{\ibf\cdot j\cdot x}
    \st j\in\Nbb_{0},\, |j|\leq k \right\rrangle
  \\
  \Tcal_{k}
  &:=
  \left\llangle \ebf^{\ibf\cdot j\cdot x} \st j\in\Nbb,\, |j|\geq k\right\rrangle
\end{align*}
where $\llangle \ldots \rrangle$ stands for the $L^{2}$-closure of the
linear span.

For any $k\in\Nbb$ the spaces $\Pcal_{k-1}$ and $\Tcal_{k}$ are
orthogonal complements of each other. Note also that
$\Tcal_{1}\subset L^{2}$ is exactly the space of functions
with zero average.

The Bernstein
inequality,~\cite{bernstein1912meilleure}, states that any function
$f\in \Pcal_{k}$ satisfies the inequality
\[
\|f'\|_{\infty}\leq k\cdot\|f\|_{\infty}
\]
with functions $\ebf^{\pm\ibf kx}$ saturating it.
Several proofs of this fact and its various generalizations can be found in
\cite{queffelec2019bernstein}. The cited article was very
inspiring while proving the \textit{reverse Bernstein inequality} below.

For $m\in\Nbb$ let $W^{m,\infty}$ be the class of functions $f\in
L^{2}(S^{1},\mu)$ such that $f\in C^{m-1}(S^{1})$ and $f$ is $m$ times
differentiable at almost every point with bounded derivative. Denote by
$\|f^{(m)}\|_{\infty}\in[0,\infty)$ the $L^{\infty}$-norm of the
  $m^{\text{th}}$ derivative of $f\in W^{m,\infty}$. For functions outside
  of this class we set $\|f^{(m)}\|_{\infty}:=\infty$.

\begin{theorem}[Reverse Bernstein Inequality]
  \label{p:rbi}\ \\
  Let $k,m$ be natural numbers and $f\in \Tcal_{k}$, then
  \[
  \|f^{(m)}\|_{\infty}\geq C_{k,m}\cdot\|f\|_{\infty}
  \]
  where $C_{k,m} = 1/\|I^{m-1} c_{k}\|_{\infty}$.
  Here $I^{m-1} c_{k}$ stands for the
  $(m-1)$-primitive of the triangular function $c_{k}$ with period
  $2\pi/k$.
  \qed
\end{theorem}

The operator $I$ and triangular functions $c_{k}$ are discussed in the next
section, see~(\ref{eq:triangular}) and~(\ref{eq:primitive-operator}).

It is easy to see, that
$C_{k,m}=\left(\frac{2k}{\pi}\right)^{m}\cdot B_{m}^{-1}$ for some
constants $B_{m}$. Closed expression for $B_{m}$ is not known to us at
the moment, however, it is easy to find the value for any given
$m\in\Nbb$. For example, $C_{k,1}=2k/\pi$ and
$C_{k,2}=8k^{2}/\pi^{2}$. The last section contains description of a
simple algorithm for evaluating $C_{k,m}$.

The inequality of Theorem~\ref{p:rbi} for $k=1$ is proven
in~\cite{northcott1939some}. In the
article~\cite{partington1983resolvent} it is shown by different
methods that for every $f\in \Tcal_{1}$ and every $p\in[1,\infty]$ the
following inequality holds
\[
  \|f'\|_{p}\geq C_{1,1}\|f\|_{p}
\]
It was communicated to us by J. Partington, that techniques of the last
cited article can be generalized to prove reverse Bernstein inequality for
every $k,m\in\Nbb$ and every Lebesgue norm.

\section{Proof of the Theorem}
We will need the following notations in the proof.
Define \textit{triangular cosine and sine} by
\begin{align}\label{eq:triangular}
  &c(x)
  :=
  \frac{\pi}2-d(x,0),
  &&
  c_{k}(x)
  :=\frac{c(k\cdot x)}{k},
  \\
  \nonumber
  &s(x)
  :=
  c(\frac\pi2-x),
  &&s_{k}(x)
  :=\frac{s(k\cdot x)}{k}
\end{align}

\long\def\triangularfig{
\begin{wrapfigure}[7]{o}{60 mm}
\centering
\begin{tikzpicture}
[yscale=0.68,xscale=0.68]
  \coordinate (O) at (0,0);
  \draw[->] (-4,0) -- (4,0) coordinate[] (xmax);
  \draw[->] (0,-1.7) -- (0,1.9) coordinate[] (ymax);
  \draw (3.14159,0) node[below]{${\pi}$};
  \draw(3.14159,0) -- (3.14159,0.1);
  \draw (-3.14159,0) node[below]{$-{\pi}$};
  \draw(-3.14159,0) -- (-3.14159,0.1);
  \draw[thick, blue] (-3.14159,-1.57079) -- (0,1.57079);
  \draw[thick, blue] (0,1.57079) -- (3.14159,-1.57079);
  \draw (3.14159,-1.57079) node[right,blue]{$c(x)$};
  \draw (1.67079,0) node[above]{${\pi}/2$};
  \draw(1.57079,0) -- (1.57079,0.1);
  \draw (-1.78,0) node[above]{$-{\pi}/2$};
  \draw(-1.57079,0) -- (-1.57079,0.1);
  \draw[thick, red] (-3.14159,0.78539) -- (-1.57079,-0.78539);
  \draw[thick, red] (-1.57079,-0.78539) -- (0,0.78539);
  \draw[thick, red] (0,0.78539) -- (1.57079,-0.78539) ;
  \draw[thick, red] (1.57079,-0.78539)  -- (3.14159,0.78539);
  \draw (3.,0.78539) node[below,right, red]{$c_2(x)$};
\end{tikzpicture}
\caption{\small Graphs~of~$c(x)$~and~$c_2(x)$.}\label{fig:triangular}
\vskip0mm
\end{wrapfigure}
}
\skipthis{
\begin{figure}
\begin{tikzpicture}
  \coordinate (O) at (0,0);
  \draw[->] (-4,0) -- (4,0) coordinate[] (xmax);
  \draw[->] (0,-1.7) -- (0,1.7) coordinate[] (ymax);
  \draw (3.14159,0) node[below]{${\pi}$};
  \draw(3.14159,0) -- (3.14159,0.1);
  \draw (-3.14159,0) node[below]{$-{\pi}$};
  \draw(-3.14159,0) -- (-3.14159,0.1);
  \draw[thick, blue] (-3.14159,-1.57079) -- (0,1.57079);
  \draw[thick, blue] (0,1.57079) -- (3.14159,-1.57079);
  \draw (3.14159,-1.57079) node[right,blue]{$c(x)$};
  \draw (1.67079,0) node[above]{${\pi}/2$};
  \draw(1.57079,0) -- (1.57079,0.1);
  \draw (-1.78,0) node[above]{$-{\pi}/2$};
  \draw(-1.57079,0) -- (-1.57079,0.1);
  \draw[thick, red] (-3.14159,0.78539) -- (-1.57079,-0.78539);
  \draw[thick, red] (-1.57079,-0.78539) -- (0,0.78539);
  \draw[thick, red] (0,0.78539) -- (1.57079,-0.78539) ;
  \draw[thick, red] (1.57079,-0.78539)  -- (3.14159,0.78539);
  \draw (3.14159,0.78539) node[below,right, red]{$c_2(x)$};
\end{tikzpicture}
\caption{Graphs of $c(x)$ and $c_2(x)$.}\label{fig:triangular}
\end{figure}}
The functions $s_k$ and $c_{k}$ are piece-wise linear and have the
same sign and monotonicity as $\sin(k\cdot x)$ and $\cos(k\cdot x)$,
respectively, and for every $k\in\Nbb$ and almost every $x\in S^{1}$
they satisfy $|c'_k(x)|=|s'_{k}(x)|=1$, see Fig~\ref{fig:triangular}.
Clearly $s_{k},c_{k}\in \Tcal_{k}$ and moreover
\[
  \Tcal_{k}= \llangle c_{i},s_{i} \st i\geq k \rrangle
\]

\triangularfig
This is because the Fourier coefficients of functions $c$ and $s$ are
multiplicative and it was shown in \cite{hartman1947multiplicative}
that in such case the dilations of the functions form a Riesz basis of
$L^{2}$, see also~\cite{wei1999triangular} for the discussion about
triangular functions specifically.

\begin{proof}[Proof of the Theorem~\ref{p:rbi}]
  Consider the bounded linear operator
  \[
    I:\Tcal_{1}\to \Tcal_{1}
  \]
  defined for $\phi\in \Tcal_{1}$ by
  \begin{equation}\label{eq:primitive-operator}
    I\phi := \Phi
  \end{equation}
  where $\Phi$ is the (unique) primitive function of $\phi$ with zero
  average.  Note that the operator $I$ leaves spaces $\Tcal_{k}$ and
  $\Tcal_{1}\cap\Pcal_{k}$ invariant for every $k\in\Nbb$, since
  the standard $L^{2}$-basis of $\Tcal_{1}$ is the eigenbasis of $I$.
    
  Define $J_{m}\in \Tcal_{1}$ by
  \begin{align*}
    J_{1}(x)
    &:=
    x
    &&\text{for $x\in(-\pi,\pi)$}
    \\
    J_{1}(\pi)&:=0
    \\
    J_{m}&:=I^{m-1}J_{1}
    &&\text{for $m>1$}
  \end{align*}
  Then for every $\phi\in\Tcal_{1}$
  \begin{align*}
    \< \phi,J_{1}\> &=
    \frac1{2\pi}\int_{-\pi}^{\pi}\phi(x)\cdot x\d x =
    \lim_{\epsilon\to0}\left.\frac1{2\pi}x\cdot\Phi(x)\right|_{-\pi+\epsilon}^{\pi-\epsilon}
    \\    \nonumber
    &=
    \Phi(\pi)
    = (I\phi)(\pi) 
  \end{align*}
  Since operator $I$ is anti-self-adjoint, we also have
  \begin{equation}\label{eq:m-primitive}
    |I^{m}\phi(\pi)|=|\< \phi,J_{m}\>|
  \end{equation}

  Now let $k,m$ be natural numbers. Take $\phi\in\Tcal_{k}$ and let
  $p\in \Pcal_{k-1}$ be arbitrary. Since $\phi$ and $p$ are orthogonal,
  Equation~(\ref{eq:m-primitive}) implies
  \[
  |(I^{m}\phi)(\pi)| = |\< \phi, (J_{m}-p) \>|
  \]

  By Hölder inequality  it
  follows that for every $p\in\Pcal_{k-1}$
  \begin{equation*}
    |(I^{m}\phi)(\pi) |\leq \|\phi\|_{\infty}\cdot \|J_{m} - p\|_{1}
  \end{equation*}

  The shift operator commutes with $I^{m}$ and applying it to
  the function $\phi$ does not affect the right-hand-side of the
  inequality. Thus we have
  \begin{equation}\label{eq:m-primitive-norm}
    \|I^{m}\phi\|_{\infty}\leq \|\phi\|_{\infty}\cdot \|J_{m} - p\|_{1}
  \end{equation}
    
  Now we need to choose $p\in\Pcal_{k-1}$ smartly. More specifically,
  we will choose $p$ in such a way that there is a resonance in
  the Hölder inequality for $\phi=c_{k}'$ or $\phi=s_{k}'$ depending
  on the parity of $m$.

  To achieve this we will take $p\in \Pcal_{k-1}$ to be the
   trigonometric Lagrange polynomial of the function $J_{m}$ at points
  uniformly distributed in the circle and make use of the
  two Lemmas~\ref{p:bounds} and~\ref{p:existence} below.

  We call a function on the circle even/odd if it is even/odd with
  respect to the complex conjugation on the circle (or involution
  $x\mapsto -x$, in our coordinates). Clearly, for any
  $m\in\Nbb$ the parity of $J_{m}$ is the same as the parity of $m$.
  Also for every $m\in\Nbb$ we have $J_{m}\in W^{m,\infty}$.
  We denote by $\Pcal_{k}^{\even}$ and  $\Pcal_{k}^{\odd}$ the space
  of even (respectively, odd) trigonometric polynomials of degree
  $k\in\Nbb_{0}$ and for notational convenience also set
  $\Pcal^{(m)}_{k}$ to be $\Pcal_{k}^{\even}$ or $\Pcal_{k}^{\odd}$,
  depending on the parity of $m$.

  \begin{lemma}\label{p:bounds}
    Let $k,m\in\Nbb$ and $p\in\Pcal_{k}^{(m)}$.
    Let
    \begin{align*}
      Z&:=\set{x\in S^{1}\st (J_{m}-p)(x)=0}\\
    \end{align*}
    Then
    \def\theenumi{\roman{enumi}}
    \begin{enumerate}
    \item\label{p:bounds:bound}
      $|Z|\leq 2k+2$
    \item\label{p:bounds:simple}
      If $|Z|= 2k+2$ then all points in $Z$ are simple zeroes of $(J_{m}-p)$.
    \end{enumerate}
    \qed
  \end{lemma}
  A remark here is in order. Function $J_{1}$ is not continuous and
  $J_{2}$ is not differentiable at $x=\pi$. In that case we say that
  $x=\pi$ \textit{is} a simple zero of $J_{1}-p$ if
  $p(\pi)=0$, and $x=\pi$ \textit{is never} a simple zero
  of $J_{2}-p$.  The cases $m=1,2$ will need a little special care in the
  proof of Lemma~\ref{p:bounds}.

  The second lemma deals with the existence and uniqueness of
  an interpolation by trigonometric polynomials. Certainly, theory of
  such interpolations is well developed, see for
  example~\cite[Chapter 3.8]{atkinson1991introduction}.
  However, we deal here with additional
  symmetries (parity of interpolated function) and since the proof is
  rather short we are compelled to include the lemma together with its
  proof here.
  
  \begin{lemma}\label{p:existence}
    Let $k,m\in\Nbb$.
    Then
    \def\theenumi{\roman{enumi}}
    \begin{enumerate}
    \item\label{p:existence:odd}
      If $m$ is odd and $A\subset(0,\pi)$ has cardinality $k$ then
      there exists unique odd polynomial $p\in\Pcal_{k}^{(m)}$ such
      that
      \[
      \set{x\in S^{1}\st (J_{m}-p)(x)=0} = A\cup(-A)\cup\set{0,\pi}
      \]
    \item\label{p:existence:even}
      If $m$ is even and $A\subset(0,\pi)$ has cardinality $k+1$ then
      there exists unique even polynomial $p\in\Pcal_{k}^{(m)}$ such
      that
      \[
      \set{x\in S^{1}\st (J_{m}-p)(x)=0} = A\cup(-A)
      \]
    \end{enumerate}
    \qed
  \end{lemma}
  
  We postpone the proofs of the lemmas until after the end of the proof of the
  theorem.

  Assume now that $m\in\Nbb$ is odd and $k\in\Nbb$. Let
  \[
    A:=\set{\frac{j\cdot\pi}{k}\st
      j=1,\dots,(k-1)}\subset(0,\pi)
  \]
  By Lemma~\ref{p:existence}(\ref{p:existence:odd}) there exists
  unique $p\in\Pcal_{k-1}^{\odd}$ iterpolating $J_{m}$ at
  $A\cup(-A)\cup\set{0,\pi}$.
  By
  Lemma~\ref{p:bounds}(\ref{p:bounds:simple}) the function $(J_{m}-p)$
  changes sign at every point of $A\cup(-A)\cup\set{0,\pi}$ and does not change sign in
  between. Since function $c'_{k}$ is unimodular and also changes sign
  at the same points, it follows that
  \begin{equation}\label{eq:l1-optimum-odd}
    \|J_{m}-p\|_{1}
    =
    |\< J_{m}-p,c'_{k}\>|
    =
    |\< J_{m}, c'_{k}\>|
    =\| I^{m}c_{k}' \|_{\infty}
    =:D_{k,m}
  \end{equation}

  For even $m$ we reason in a similar way. Set
  \[
    A:=\set{\frac{(j+1/2)\cdot\pi}{k}\st
      j=0,\dots,k-1}\subset(0,\pi)
  \]
  and also let $p\in\Pcal_{k}^{\even}$ be the trigonometric
  polynomial provided by
  Lemma~\ref{p:existence}(\ref{p:existence:even}).  Then by
  Lemma~\ref{p:bounds}(\ref{p:bounds:simple}) functions $J_{m}-p$
  and $s'_{k}$ change sign at exactly the same points and we have
  \begin{equation}\label{eq:l1-optimum-even}
    \|J_{m}-p\|_{1}
    =
    |\< J_{m}-p,s'_{k}\>|
    =
    |\< J_{m}, s'_{k}\>|
    =\| I^{m}s'_{k} \|_{\infty}
    =:D_{k,m}
  \end{equation}
 
  Since for $\phi=c'_{k}$ or $\phi=s'_{k}$ depending on the parity of
  $m$ we have equality in the inequality~(\ref{eq:m-primitive-norm}),
  we have also proven here that such interpolating polynomial $p$ is
  $L^{1}$-optimal.
  
  Combining the inequality~(\ref{eq:m-primitive-norm})
  with~(\ref{eq:l1-optimum-odd}) or~(\ref{eq:l1-optimum-even})
  we obtain the inequality
  \begin{equation}\label{eq:rbi-integral}
  \| I^{m}\phi\|_{\infty}
  \leq
  D_{k,m}\cdot\|\phi\|_{\infty}
  \end{equation}

  To finish the proof of the theorem, let $k,m\in\Nbb$ and $f\in\Tcal_{k}$.
  If $f\not\in W^{m,\infty}$, the conclusion of the theorem is trivially
  satisfied. For $f\in W^{m,\infty}$ we have $f=I^{m}f^{(m)}$.
  Applying the inequality~(\ref{eq:rbi-integral}) to the function
  $\phi:=f^{(m)}$ we obtain the conclusion of the theorem with
  $C_{k,m}=\left(D_{k,m}\right)^{-1}$. 
  It remains only to note that the function $I^{m}c'_{k}$ saturates the
  inequality and the constant $C_{k,m}$ is optimal. 
\end{proof}

\begin{proof}[Proof of Lemma~\ref{p:bounds}]  
  First we observe that the set $Z':=Z\cap(-\pi,\pi)$ is a discrete subset of
  $(-\pi,\pi)$, since $(J_{m}-p)$ is analytic on $(-\pi,\pi)$. Also
  $Z$ is invariant under the involution $x\mapsto -x$ on the
  circle. If $m$ is odd then $\set{0,\pi}\subset Z$.

  To prove the assertion~(\ref{p:bounds:bound}) of the lemma we
  proceed by induction with respect to $m$.

  \underline{Let $m=1$} and $p\in\Pcal_{k}^{\odd}$.  The function
  $J_{1}-p$ has at least one critical point strictly between any
  two consecutive points in $Z'$.  There are at most $2k$ critical
  points of $J_{1} - p$, since $(J_{1}-p)'=1-p'$ is a non-zero
  trigonometric polynomial of degree not greater then $k$ and may not
  have more then $2k$ zeroes on $S^1$. It follows that $|Z'|\leq 2k+1$
  and therefore $|Z|\leq 2k+2$.

  \underline{Let $m>1$}, then $(J_{m}-p)$ is continuous on the circle
  and analytic on $(-\pi,\pi)$. It has a critical point (possibly
  including a point $x=\pi$ where derivative is not defined for $m=2$)
  strictly between any two consecutive points in $Z$.  But by
  induction assumption there are at most $2k+2$ zeroes of
  $(J_{m-1}-p')=(J_{m}-p)'$, therefore $|Z|\leq 2k+2$.

  To prove the assertion~(\ref{p:bounds:simple}) observe that under
  the assumptions of the lemma the function $(J_{m}-p)$ has exactly
  $2k+2$ zeroes and at most $2k+2$ critical points. Since there must
  be a critical point strictly between any two consecutive zeroes,
  there are no critical points to spare and all zeroes must be simple.
\end{proof}

\begin{proof}[Proof of Lemma~\ref{p:existence}]
First we note that independently of the parity of $m$ any two
interpolating polynomials in $\Pcal_{k}$ will have the same values at
$2k+2$ distinct points, therefore their difference will be identically
zero. Thus, interpolating polynomial is unique, if exists.

To prove existence we treat cases of even and odd $m$ separately.
Let $m$ be even. Define for every $a\in A$
\[ 
p_{a}(x):= \prod_{b\in A\setminus\set{a}}\big(\cos(x)-\cos(b)\big)
\]
Since for each $a\in A$ polynomial $p_{a}$ is even and $p_{a}(a)\neq0$,
we can find interpolating polynomial as
\[
p(x):=\sum_{a\in A}J_{m}(a)\frac{p_{a}(x)}{p_{a}(a)}
\]
For odd $m$ consider the set $Z':=A\cup(-A)\cup\set{\pi}$
and define for every $z\in Z'$
\[
p_{z}(x)
:=
\ebf^{-\ibf\cdot k\cdot x}
\prod_{\zeta\in Z'\setminus\set{z}}(\ebf^{\ibf\cdot x}-\ebf^{\ibf\cdot \zeta})
\]
Clearly, $p_{z}(z)\neq0$ and we define the interpolating polynomial
\[
p(x):=\sum_{z\in Z'}J_{m}(z)\frac{p_{z}(x)}{p_{z}(z)}
\]
It will have the required values at all points in $Z'$. Since $-p(-x)$
and $\overline{p(x)}$
will also have the same values at points in $Z'$ it follows that $p$
is odd, real polynomial and therefore
\begin{align*}
  p(0)&=J_{m}(0)=0\\
  p(\pi)&=J_{m}(\pi)=0
\end{align*}
Thus $p$ interpolates $J_{m}$ at $A\cup(-A)\cup\set{0,\pi}$.
\end{proof} 

\section{The constants}
Here we give a simple algorithm for evaluating the constants
\[
  C_{k,m}=\frac{1}{\|I^{m}c'_{k}\|_{\infty}}
\]
for $k,m\in\Nbb$. The derivation of the algorithm is elementary and
not very enlightening, thus it is not included.

Define a sequence of algebraic polynomials on $[0,1]$ recursively by
\begin{align*}
  P_{0}(x)&:=1\\
  P_{m+1}(x)&:=
            \begin{cases}
              \int_{0}^{x}P_{m}\drm x & \text{$m$ is even}\\
              \int_{x}^{1}P_{m}\drm x & \text{$m$ is odd}
            \end{cases}           
\end{align*}

Define constants $B_{m}$ as
\begin{align*}
  B_{m}:=
  \begin{cases}
    P_{m}(0) & \text{$m$ is even}\\
    P_{m}(1) & \text{$m$ is odd}
  \end{cases}
\end{align*}

At least 27 initial terms of the
sequence $(B_{m}\cdot m!)_{0}^{\infty}$ match the sequence $A000111$,
\textit{Euler or up/down numbers}, in \url{https://oeis.org/A000111}.

\begin{proposition}
  For $k,m\in\Nbb$
  \[
    C_{k,m}:=\left(\frac{2k}{\pi}\right)^{m} B_{m}^{-1}
  \]
  In particular,
  \begin{align*}
    C_{k,1}=\frac{2k}{\pi}
      \quad\text{and}\quad
    C_{k,2}=\frac{8k^{2}}{\pi^{2}}.
  \end{align*}
  \qed
\end{proposition}


\printbibliography

\end{document}